# Optimal portfolio selection and compression in an incomplete market [*]

Nikolai Dokuchaev[†] and Ulrich Haussmann [‡]

**Abstract**

We investigate an optimal investment problem with a general performance criterion which, in particular, includes discontinuous functions. Prices are modeled as diffusions and the market is incomplete. We find an explicit solution for the case of limited diversification of the portfolio, i.e. for the portfolio compression problem. By this we mean that an admissible strategies may include no more than $m$ different stocks concurrently, where $m$ may be less than the total number $n$ of available stocks.

## 1 Introduction

The paper investigates an optimal investment problem for a market which consists of a risk free bond and a finite number of risky stocks. It is assumed that the stock prices evolve according to an Itô stochastic differential equation. The problem goes back to Black and Scholes (1973) and Merton (1969,1973). Black and Scholes gave strategies to replicate a given claim, while Merton found the strategies which solve an optimization problem in which $\mathbf{E}U(X(T))$ is to be maximized, where $X(T)$ represents the wealth at final time $T$ and $U(\cdot)$ is a utility function (see e.g. Laurent and Pham (1999), Khanna and Kulldorff (1999), Zhou (1998)).

[*]Supported by NSERC under NCE Grant 30354 and Research Grant 88051, and by Russian Foundation for Basic Research grant 99-01-0886.

[†]The Institute of Mathematics and Mechanics, St.Petersburg State University; since 2001, Department of Mathematics and Computer Science, The University of West Indies, Mona, Jamaica. Email: ndokuch@uwimona.edu.jm; N_Dokuchaev@yahoo.com

[‡]Department of Mathematics, University of British Columbia, 1984 Mathematics Road, Vancouver, BC, Canada, V6T 1Z2. Email: uhaus@math.ubc.ca

In this paper we study the investment problem for a general type utility function which covers utilities of quadratic form, log form, power form and discontinuous utilities. But we also consider the case of limited diversification of the portfolio, i.e. the portfolio compression problem. By this we mean that an admissible strategies may include no more than $m$ different stocks concurrently, where $m$ may be less than the total number $n$ of available stocks. Although this problem has not been treated in the literature, it is of interest to the investor. Clearly, it is not realistic to include in his portfolio *all* available stocks; the total number of assets in the market is too large. In fact, the number of stocks in the portfolio should be limited by the equity in the account (say, several hundred stocks for a large fund, and fewer for an individual investor), because of the need to have a position in each stock large enough so that management fees and commissions are only a small proportion of the value of the portfolio. There is just not much point in having too many stocks in a small portfolio. Even in a large portfolio it makes sense to limit the number of stocks to those that can be watched closely. On the other hand, there should be a certain minimum number of stocks so that a sufficient degree of diversification can be achieved. The criteria for making the selection will depend on the investor's preferences; we will assume that a number $m$ is fixed a priori and then the $m$ stocks that yield the best return are sought. We obtain the optimal strategy explicitly under some conditions.

Our model consists of a multi-stock diffusion market with correlated stock prices. The parameters (interest rate, appreciation rate and volatilities) need not be adapted to the driving Wiener process, so the market is incomplete. We treat three cases: (i) the utility function is the log, (ii) the market parameters are independent of the driving Wiener process (but still random) and the $L_2$ norm (in time) of the market price of risk is nonrandom, and (iii) the market parameters are independent of the driving Wiener process and the utility function is of special type. For the case $n = m$, this last result extends a result of Karatzas and Shreve (1998) on "totally unhedgable" coefficients. We also exhibit the optimal portfolio explicitly in terms of the state-price density.

In Section two we collect notation and definitions, and we set up the market. The problem is stated in Section three, and in Section four a formula for replicating some special claims is presented. In Section five, the problem is solved when all the stocks maybe held in the portfolio, and in Section six it is solved when not all stocks may be held. The proofs are given in the Appendix.

## 2    The Model and Definitions

Consider a diffusion model of a market consisting of a risk free bond or bank account with the price $B(t)$, $t \geq 0$, and $n$ risky stocks with prices $S_i(t)$, $t \geq 0$, $i = 1, 2, ..., n$, where $n < +\infty$ is given. The prices of the stocks evolve according to the following equations:

$$dS_i(t) = S_i(t) \left( a_i(t)dt + \sum_{j=1}^{n} \sigma_{ij}(t)dw_j(t) \right), \quad t > 0, \tag{2.1}$$

where the $w_i(t)$ are standard independent Wiener processes, $a_i(t)$ are appreciation rates, and $\sigma_{ij}(t)$ are volatility coefficients. The initial price $S_i(0) > 0$ is a given non-random constant. The price of the bond evolves according to the following equation

$$B(t) = B(0) \exp \left( \int_0^t r(s)ds \right), \tag{2.2}$$

where $B(0)$ is a given constant which we take to be 1 without loss of generality, and $r(t)$ is a random process of risk-free interest rate.

We are given a standard probability space $(\Omega, \mathcal{F}, \mathbf{P})$, where $\Omega$ is the set of all events, $\mathcal{F}$ is a complete $\sigma$-algebra of events, and $\mathbf{P}$ is a probability measure. Introduce the vector processes ($^\top$ denoted transpose)

$$w(t) = (w_1(t), ..., w_n(t))^\top, \quad a(t) = (a_1(t), ..., a_n(t))^\top, \quad S(t) = (S_1(t), ..., S_n(t))^\top$$

and the matrix process $\sigma(t) = \{\sigma_{ij}(t)\}_{i,j=1}^n$.

We assume that $\{w(t)\}_{0 \leq t \leq T}$ is a standard Wiener process, $a(t)$, $r(t)$ and $\sigma(t)$ are uniformly bounded, measurable random process, independent of future increments of $w$, such that $c_1 \mathbf{I}_n \leq \sigma(t)\sigma(t)^\top$ where $c_1 > 0$ is a constant and $\mathbf{I}_n$ is the identity matrix in $\mathbf{R}^{n \times n}$. Under these assumptions the solution of (2.1) is well defined, but the market is incomplete.

Set
$$V(t) \triangleq \sigma(t)\sigma(t)^\top, \quad Q(t) \triangleq V(t)^{-1}$$

$$\widehat{r}(t) \triangleq (r(t), r(t), ..., r(t))^\top, \quad \widetilde{a}(t) \triangleq a(t) - \widehat{r}(t), \quad \theta(t) \triangleq \sigma(t)^{-1}\widetilde{a}(t),$$

$$R \triangleq \int_0^T |\theta(t)|^2 dt, \quad \overline{R} \triangleq \frac{R}{T}, \quad J \triangleq \mathbf{E}R, \quad \tau(t) \triangleq \overline{R}^{-1} \int_0^t |\theta(s)|^2 \, ds, \tag{2.3}$$

and let $\mathbf{S}(t) \triangleq \mathrm{diag}\,(S_1(t), ..., S_n(t))$ be the diagonal matrix with the corresponding diagonal elements.

Let $\mu(t) \triangleq (r(t), a(t), \sigma(t))$. Let $\{\mathcal{F}_t\}_{0 \leq t \leq T}$ be the filtration generated by the process $(S(t), \mu(t))$ completed with the null sets of $\mathcal{F}$. By (2.1),

$$dw(t) = \sigma(t)^{-1}\mathbf{S}(t)^{-1} \left(dS(t) - \mathbf{S}(t)a(t)dt\right). \tag{2.4}$$

It follows that $\{\mathcal{F}_t\}$ coincides with the filtration generated by the processes $(w(t), \mu(t))$.

Set
$$p(t) \triangleq \exp\left(-\int_0^t r(s)ds\right) = B(t)^{-1}, \quad \widetilde{S}(t) \triangleq p(t)S(t). \tag{2.5}$$

It is easy to see that $\mathcal{F}_t$ coincides with the filtration generated by the processes $(\widetilde{S}(t), \mu(t))$.

Let
$$w_*(t) \triangleq w(t) + \int_0^t \theta(s)ds,$$

and
$$\mathcal{Z}_*(t) = \exp\left(-\int_0^t \theta(s)^\top dw(s) - \frac{1}{2}\int_0^t |\theta(s))|^2 ds\right). \tag{2.6}$$

Our standing assumptions imply that $\mathbf{E}\mathcal{Z}_*(T) = 1$. Define the (equivalent martingale) probability measure $\mathbf{P}_*$ by $d\mathbf{P}_*/d\mathbf{P} = \mathcal{Z}_*(T)$. Let $\mathbf{E}_*$ be the corresponding expectation. Girsanov's theorem implies that $w_*$ is a standard Wiener process under $\mathbf{P}_*$. Then $(w_*(t), \mu(t))$ also generate $\{\mathcal{F}_t\}$ and $\mathbf{P}$, $\mathbf{P}_*$ have the same null sets. Now define

$$\mathcal{Z}(t) = \exp\left(\int_0^t \theta(s)^\top dw_*(s) - \frac{1}{2}\int_0^t |\theta(s))|^2 ds\right). \tag{2.7}$$

It follows from (2.7) that

$$dZ(t) = Z(t)\theta(t)^\top dw_*(t) = Z(t)\widetilde{a}(t)^\top Q(t)\widetilde{\mathbf{S}}(t)^{-1}d\widetilde{S}(t), \qquad (2.8)$$

where $\widetilde{\mathbf{S}}(t) \triangleq \operatorname{diag}(\widetilde{S}_1(t), ..., \widetilde{S}_n(t))$. But $Z_*(t) = Z(t)^{-1}$ since Itô's Lemma implies that $d(Z_*(t)Z(t)) = 0$. Note that $\mathbf{E}_*Z(T) = 1$, $2\mathbf{E}_* \log Z(T) = -J$.

Let $X_0 > 0$ be the initial wealth at time $t = 0$, and let $X(t)$ be the wealth at time $t > 0$, $X(0) = X_0$. We assume that

$$X(t) = \pi_0(t) + \sum_{i=1}^n \pi_i(t), \qquad (2.9)$$

where the pair $(\pi_0(t), \pi(t))$ describes the portfolio at time $t$. The process $\pi_0(t)$ is the investment in the bond, $\pi_i(t)$ is the investment in the $i$th stock, $\pi(t) = (\pi_1(t), ...., \pi_n(t))^\top$, $t \geq 0$.

The portfolio is said to be self-financing, if

$$dX(t) = \pi(t)^\top \mathbf{S}(t)^{-1} dS(t) + \pi_0(t) B(t)^{-1} dB(t). \qquad (2.10)$$

It follows that for such portfolios

$$dX(t) = r(t)X(t)\,dt + \pi(t)^\top \left(\widetilde{a}(t)\,dt + \sigma(t)\,dw(t)\right), \qquad (2.11)$$

$$\pi_0(t) = X(t) - \sum_{i=1}^n \pi_i(t),$$

so $\pi$ alone suffices to specify the portfolio; it is called a self-financing strategy.

**Definition 2.1** *The process $\widetilde{X}(t) \triangleq p(t)X(t)$ is called the normalized wealth.*

It satisfies

$$\widetilde{X}(t) = X(0) + \int_0^t p(s)\pi(s)^\top \sigma(s)\,dw_*(s). \qquad (2.12)$$

**Definition 2.2** *Let $\overline{\Sigma}$ be the class of all $\mathcal{F}_t$-adapted processes $\pi(\cdot)$ such that for a sequence of stopping times, $\{T_k\}$ with $T_k \uparrow T$ a.s.*

- $\int_0^{T_k} \left(|\pi(t)^\top \widetilde{a}(t)| + |\pi(t)^\top \sigma(t)|^2\right) dt < \infty$ *a.s.*

- $X(T) \triangleq \lim_{k\to\infty} X(T_k)$ exists a.s.

- $\mathbf{E}_* \widetilde{X}(T) = X(0)$.

A process $\pi(\cdot) \in \overline{\Sigma}$ is said to be an *admissible* strategy with corresponding wealth $X(\cdot)$. Of course if the first condition in Definition 2.2 holds with $T_k = T$, then the other two are redundant. It turns out that the replicating strategies we use are given by the first spatial derivative of the solution of the heat equation, so may not be sufficiently regular at $t = T$ to allow us to take $T_k = T$. For an admissible strategy $\pi(\cdot)$, $X(t, \pi(\cdot))$ denotes the corresponding total wealth, and $\widetilde{X}(t, \pi(\cdot))$ the corresponding normalized total wealth.

The following definition is standard.

**Definition 2.3** *Let $\xi$ be a given random variable. An admissible strategy $\pi(\cdot) \in \overline{\Sigma}$ is said to replicate the claim $\xi$ if $X(T, \pi(\cdot)) = \xi$    a.s.*

Now we define a class of strategies where the portfolios at any one time may contain no more than a predetermined number of securities. Let $m$ be a given integer, $1 \leq m \leq n$ and let $\mathcal{M}_m$ be the collection of subsets of $\{1, ...., n\}$ each of which contains at most $m$ elements. Let $\widetilde{\mathcal{I}}_m$ be the set of $\mathcal{F}_t$-adapted functions $I : [0, T] \times \Omega \mapsto \mathcal{M}_m$, and let $\mathcal{I}_m$ be a given subset of $\widetilde{\mathcal{I}}_m$.

**Definition 2.4** *Let $\overline{\Sigma}(m)$ be the set of all strategies $\pi(\cdot) \in \overline{\Sigma}$ such that there exists $I \in \mathcal{I}_m$ for which $\pi_i(t) = 0$ if $i \notin I(t)$.*

For a given $I(\cdot) \in \mathcal{I}_m$ and $t \in [0, T]$, we denote by $L_I(t)$ the linear subspace of $\mathbf{R}^n$ such that $x = (x_1, ..., x_n) \in L_I(t)$ if and only if $x_i = 0$ for all $i \notin I(t)$.

Let $P_I(t) \in \mathbf{R}^{n\times n}$ be the projection of $\mathbf{R}^n$ onto $L_I(t)$. In other words, $P_I(t) = \{P_I^{(i,j)}(t)\}_{i,j=1}^n$, where
$$P_I^{(i,j)}(t) \triangleq \begin{cases} 1 & \text{if } i=j,\ i \in I(t) \\ 0 & \text{if } i \neq j \text{ or } i \notin I(t). \end{cases}$$

It follows that $b^\top P_I(t)V(t)P_I(t)b \geq c_1|P_I(t)b|^2$ ($\forall b \in \mathbf{R}^n$), so $P_I(t)V(t)P_I(t)$ is invertible on $L_I(t)$. Hence there exists a unique matrix $Q_I(t) = \{Q_I^{(i,j)}(t)\}_{i,j=1}^n \in \mathbf{R}^{n\times n}$ such that $Q_I^{(i,j)}(t) = 0$ if $i \notin I(t)$ or $j \notin I(t)$ and $P_I(t)V(t)P_I(t)Q_I(t)x = P_I(t)V(t)Q_I(t)x = x$ for all $x \in L_I(t)$. Set
$$a_I(t) \triangleq \widehat{r}(t) + V(t)Q_I(t)P_I(t)(a(t) - \widehat{r}(t)). \tag{2.13}$$

It is easy to see that

$$P_I(t)a_I(t) = P_I(t)a(t), \quad P_I(t)\tilde{a}_I(t) = P_I(t)\tilde{a}(t). \tag{2.14}$$

The following results are obvious.

**Lemma 2.1** *If $\sigma_{ij}(t)$ is diagonal, then*

$$a_I(t) = \begin{cases} a_i(t) & \text{if } i \in I(t) \\ r(t) & \text{if } i \notin I(t). \end{cases}$$

**Lemma 2.2** *If $m = 1$, then*

$$Q_I^{(i,j)}(t) = \begin{cases} \left(\sum_{k=1}^n \sigma_{ik}(t)^2\right)^{-1} & \text{if } i = j \in I(t), \\ 0 & \text{otherwise}. \end{cases}$$

In fact, suppose $I(t) = \{1, \ldots, m\}$, then

$$P_I(t) = \begin{pmatrix} \mathbf{I}_m & 0 \\ 0 & 0 \end{pmatrix}, \quad V(t) = \begin{pmatrix} V_1(t) & V_3(t) \\ V_3^\top(t) & V_2(t) \end{pmatrix},$$

and

$$Q_I(t) = \begin{pmatrix} V_1^{-1}(t) & 0 \\ 0 & 0 \end{pmatrix}.$$

We shall use the notation $\tilde{a}_I$, $R_I$, $\overline{R}_I$, $J_I$, $\tau_I$, $\mathbf{P}_{*I}$, and $\mathcal{Z}_I$ for $\tilde{a}$, $\overline{R}$, $R$, $J$, $\tau$, $\mathbf{P}_*$ and $\mathcal{Z}$ defined for the corresponding $a(\cdot) = a_I(\cdot)$, but with $\sigma$ unchanged.

## 3  Problem statement

Let $T > 0$ and let $m > 0$ be an integer. Let $\widehat{D} \subset \mathbf{R}$ and $X_0 \in \widehat{D}$ be given. Let $U(\cdot) : \widehat{D} \to \mathbf{R} \cup \{-\infty\}$ such that $U(X_0) > -\infty$.

We may state our general problem as follows: Find an admissible self-financing strategy $\pi(\cdot)$ which solves the following optimization problem:

$$\text{Maximize} \quad \mathbf{E}U(\widetilde{X}(T, \pi(\cdot))) \quad \text{over} \quad \pi(\cdot) \in \overline{\Sigma}(m) \tag{3.1}$$

$$\text{subject to} \quad \begin{cases} \widetilde{X}(0, \pi(\cdot)) = X_0, \\ \widetilde{X}(T, \pi(\cdot)) \in \widehat{D} \quad \text{a.s.} \end{cases} \tag{3.2}$$

The following assumptions suffice to solve the problem when $m = n$.

**Condition 3.1** *There exists a measurable set $\Lambda \subseteq \mathbf{R}$, and a measurable function $F(\cdot, \cdot) : (0, \infty) \times \Lambda \to \widehat{D}$ such that for each $z > 0$, $\widehat{x} = F(z, \lambda)$ is a solution of the optimization problem*

$$\text{Maximize} \quad zU(x) - \lambda x \quad \text{over } x \in \widehat{D}. \tag{3.3}$$

*Moreover, this solution is unique for a.e. $z > 0$.*

**Condition 3.2** *There exist $\lambda_J \in \Lambda$, $C > 0$ and $c_0 \in (0, 1/(2J))$, such that $F(\cdot, \lambda_J)$ is piecewise continuous on $(0, \infty)$, $F(\mathcal{Z}(T), \lambda_J)$ is $\mathbf{P}_*$-integrable, and*

$$\begin{cases} \mathbf{E}_*\{F(\mathcal{Z}(T), \lambda_J)\} = X_0, \\ |F(z, \lambda_J)| \leq C z^{c_0 \log z} \quad \forall z > 0. \end{cases} \tag{3.4}$$

**Condition 3.3** *At least one of the following conditions holds:*

(i) *$U(x) \equiv \log(x)$ and $\widehat{D} = [0, +\infty)$; or*

(ii) *The random variable $R$ is constant.*

**Condition 3.4** *$F(x, \lambda) = C_1 \left(\frac{x}{\lambda}\right)^\nu + C_0$, where $C_1 \neq 0$, $C_0$ and $\nu \neq 0$ are constants, and the random variable $R$ and the process $w(\cdot)$ are independent.*

*Remark 3.1.* It is clear that Condition 3.1 is required to allow maximization of the Lagrangian. Condition 3.2 ensures that the optimal terminal wealth is replicable and Condition 3.3 allows us to find the optimal, i.e. replicating, strategy explicitly. Condition 3.4 is a useful weakening of Condition 3.3(ii) for special utility functions.

It is easy to see that these conditions are satisfied in many examples.

**Lemma 3.1** *Conditions 3.1 and 3.2 hold in the following cases:*

(i) *(logarithmic utility) $\widehat{D} = [0, +\infty)$, $U(x) = \ln(x)$, $X_0 > 0$, $\Lambda = (0, \infty)$, $F(z, \lambda) = z/\lambda$, $\lambda_J = 1/X_0$.*

(ii) (power utility) Assume $R$ is constant (so $R = J$). $\widehat{D} = [0, +\infty)$, $U(x) = \frac{1}{\delta}x^\delta$, $X_0 > 0$, where $\delta < 1$, $\delta \neq 0$, $\Lambda = (0, \infty)$, $F(z, \lambda) = (z/\lambda)^{\frac{1}{1-\delta}}$, $\lambda_J = X_0^{\delta-1} \exp\{\frac{\delta}{1-\delta} \frac{R}{2}\}$.

(iii) (mean-variance utility) Assume $R$ is constant. $\widehat{D} = \mathbf{R}$, $U(x) = -kx^2 + cx$, where $k > 0$, $c \geq 0$, $\Lambda = \mathbf{R}$, $F(z, \lambda) = (c - \lambda/z)/(2k)$, $\lambda_J = (c - 2kX_0)e^{-R}$.

(iv) Assume $R$ is constant. $\widehat{D} = [0, +\infty)$, $U(x) = -x^\delta + x$, where $\delta = 1 + 1/l$, and $l > 0$ is an integer, $X_0 > \delta^{-l}$, $\Lambda = (-\infty, 0]$, $F(z, \lambda) = (1 - \lambda/z)^l \delta^{-l}$, $\lambda_J$ is a (negative) zero of a polynomial of degree $l$.

(v) (goal achieving utility) Assume $R$ is constant. $\widehat{D} = [0, \infty)$ and

$$U(x) = \begin{cases} 0 & \text{if } 0 \leq x < \alpha, \\ 1 & \text{if } x \geq \alpha, \end{cases}$$

$0 < X_0 < \alpha$, $R > 0$, $\Lambda = (0, \infty)$ and

$$F(z, \lambda) = \begin{cases} = \alpha & \text{if } 0 < \lambda < z/\alpha, \\ \in \{0, \alpha\} & \text{if } \lambda = z/\alpha, \\ = 0 & \text{if } \lambda > z/\alpha, \end{cases}$$

$\lambda_J$ is the solution of ($\Phi$ is the cumulative of the normal distribution)

$$\Phi\left(\frac{\log \lambda}{\sqrt{R}} + \frac{\log \alpha}{\sqrt{R}} + \frac{\sqrt{R}}{2}\right) = 1 - \frac{X_0}{\alpha}.$$

For case (v) Condition 3.2 fails if $R = 0$. Note that the boundedness of the coefficients $\mu$, implies that $R$ is bounded and hence $\mathbf{E}_* \mathcal{Z}(T)^q < \infty$ for any $q \in \mathbf{R}$. This is sufficient for the integrability of $F(\mathcal{Z}(T), \lambda_J)$ in the above cases. Note also that Condition 3.4 is satisfied in cases (i)-(iii).

## 4  Replicating special claims

We show that claims of a certain kind (the kind that we need for our optimization problem) are among those that can be replicated in our incomplete market. Moreover, the assumption that the $L^2$-norm (in time) of the market price of risk is non-random allows us to exhibit the replicating strategy explicitly using a transformation of the heat equation.

**Lemma 4.1** *Let $R > 0$ be constant and let $f(\cdot) : (0, \infty) \to \mathbf{R}$ be a piecewise continuous function such that $|f(x)| \leq Cx^{c_0 \log x}$ ($\forall x > 0$), where $C > 0$ and $c_0 \in (0, (2J)^{-1})$ are constants.*

*(a) The Cauchy problem*

$$\begin{cases} \frac{\partial H}{\partial t}(x,t) + \frac{\overline{R}}{2}x^2 \frac{\partial^2 H}{\partial x^2}(x,t) = 0, \\ H(x,T) = f(x) \end{cases} \quad (4.1)$$

*has a unique solution a solution $H(\cdot) \in C^{2,1}((0,\infty) \times (0,T))$, with $H(x,t) \to f(x)$ a.e. as $t \uparrow T$.*

*(b) If in addition $f(\mathcal{Z}(T))$ is $\mathbf{P}_*$-integrable, then there exists a self-financing admissible strategy $\pi(\cdot) \in \overline{\Sigma}$, with corresponding wealth $X(t)$, which replicates the claim $B(T)f(\mathcal{Z}(T))$. $\pi$ and $X$ are given by*

$$\pi(t)^\top = B(t)\frac{\partial H}{\partial x}(\mathcal{Z}(t), \tau(t))\mathcal{Z}(t)\widetilde{a}(t)^\top Q(t), \quad X(t) = B(t)H(\mathcal{Z}(t), \tau(t)). \quad (4.2)$$

*where the function $H(\cdot, \cdot) : (0, \infty) \times [0, T] \to \mathbf{R}$ is the solution of (4.1). Moreover*

$$X(0) = \mathbf{E}_* f(\mathcal{Z}(T)), \quad (4.3)$$

Note that $C^{2,1}((0,\infty) \times (0,T))$ denotes the set of functions defined on $(0,\infty) \times (0,T)$ which are continuous and have two continuous derivatives in the first variable and one in the second.

## 5 Optimal strategy for $m = n$

We solve here the optimal portfolio selection problem, but in our incomplete market. We find that in our setting there is no hedging of the coefficients. Let us explain. In the setting generally assumed in finance, cf. Merton (1990), Sec. 15.5, the coefficients, $\mu = (r, a, \sigma)$, are assumed to satisfy an Itô equation with driving Brownian motion $(w(\cdot), \widehat{w}(\cdot))$, i.e.

$$d\mu(t) = \beta(B(t), S(t), \mu(t), t)dt + \sigma^{\mu,S}(B(t), S(t), \mu(t), t)dw(t) + \sigma^\mu(B(t), S(t), \mu(t), t)d\widehat{w}(t).$$

To Markovianize the problem, we use the state variables $X(t), B(t), S(t), \mu(t)$. Then the Bellman equation, satisfied formally by the value function (derived utility function), $J(x, s, \mu, t)$, is (we denote the matrix diag($s$) by $\mathbf{s}$)

$$\begin{aligned}
\max_{\pi} \{ & J_t(x,b,s,\mu,t) + J_x(x,b,s,\mu,t)[rx + \pi^\top(a - \hat{r})] + J_b(x,b,s,\mu,t)rb \\
& + J_s(x,b,s,\mu,t)^\top \mathbf{s}a + J_\mu(x,b,s,\mu,t)^\top \beta(b,s,\mu,t) \\
& + \tfrac{1}{2} J_{x,x}(x,b,s,\mu,t)\pi^\top \sigma\sigma^\top \pi + \tfrac{1}{2}\mathrm{tr}\,[J_{s,s}(x,b,s,\mu,t)\mathbf{s}\sigma\sigma^\top \mathbf{s}] \\
& + \tfrac{1}{2}\mathrm{tr}\,[J_{\mu,\mu}(x,b,s,\mu,t)(\sigma^{\mu,S}(b,s,\mu,t)\sigma^{\mu,S}(b,s,\mu,t)^\top \\
& \qquad + \sigma^\mu(b,s,\mu,t)\sigma^\mu(b,s,\mu,t)^\top)] + J_{x,s}(x,b,s,\mu,t)\mathbf{s}\sigma\sigma^\top \pi \\
& + J_{x,\mu}(x,b,s,\mu,t)\sigma^{\mu,S}\sigma^\top \pi + \mathrm{tr}\,[J_{s,\mu}(x,b,s,\mu,t)\sigma^{\mu,S}(b,s,\mu,t)\sigma^\top \mathbf{s}] \} = 0,
\end{aligned}$$

$$J(x, b, s, \mu, T) = U(x/b).$$

Then the optimal $\pi$ is (formally)

$$\begin{aligned}
\pi(t) = & -\frac{J_x(X(t), B(t), S(t), \mu(t), t)}{J_{x,x}(X(t), B(t), S(t), \mu(t), t)} Q(t)\widetilde{a}(t) - \mathbf{S}(t)\frac{J_{s,x}(X(t), B(t), S(t), \mu(t), t)}{J_{x,x}(X(t), B(t), S(t), \mu(t), t)} \\
& - Q(t)\sigma(t)\sigma^{\mu,S}(B(t), S(t), \mu(t), t)^\top \frac{J_{\mu,x}(X(t), B(t), S(t), \mu(t), t)}{J_{x,x}(X(t), B(t), S(t), \mu(t), t)}.
\end{aligned}$$

The first term on the right side gives the usual mean-variance type of strategy, the second, due to correlation between wealth and stock prices, is absent if $S(t)$ is not required as a state variable, e.g. if a Mutual Fund theorem holds, and the third term depends on the correlation between $S$ (or $w$) and $\mu$ and is considered to represent a hedge against future unfavourable behaviour of the coefficients. Note that the Bellman equation is degenerate: the coefficient matrix for the second order derivatives has rank at most $1 + 2n + n^2$ whereas there are $3 + 2n + n^2$ variables. The difference of 2 in the numbers arises from including $B(t)$ as a state variable (this might be avoided in some cases), and from the fact that the noise driving $X(t)$ is the same as that driving $S(t)$. This is unavoidable. Hence there may not exist a solution $J$ with second-order derivatives. If $\mu$ is independent of $w$, then $\sigma^{\mu,S} = 0$ and $B, S$ can be dropped as state variables. In this case the coefficients are said to be unhedgable and the policy "myopic". Our setting is more general in that the coefficients are not necessarily generated by a diffusion. We show rigorously that in the iso-elastic utility case, i.e. when $U(x) = x^\delta$ ($\delta < 1$) or $U(x) = \log x$, then the correlation terms are absent *provided* (in the power case) that the random variable $R$ is independent of $w$. We do not require that the coefficients $\mu$ be independent of $w$. Moreover, *if $R$ is constant*, then again these terms are absent without specifying the form of $U$.

Let $U^+(x) \triangleq \max(0, U(x))$ and let $F(\cdot, \cdot)$ be as in Condition 3.1.

**Theorem 5.1** (i) Let $R = 0$, then the trivial strategy, $\pi(t) \equiv 0$, is the unique optimal strategy for the problem (3.1)-(3.2).

(ii) Assume Conditions 3.1, 3.2 and 3.3, and let $R > 0$.

(a) If Condition 3.3(i) holds, let $H(x,t) \triangleq x/\lambda_J$;

(b) if Condition 3.3(ii) holds, let $H(x,t)$ be the solution of (4.1) for $f(\cdot) = F(\cdot, \lambda_J)$.

Then the strategy $\pi(\cdot) \in \overline{\Sigma}$ defined as

$$\pi(t)^\top = B(t)\frac{\partial H}{\partial x}(\mathcal{Z}(t), \tau(t))\mathcal{Z}(t)\widetilde{a}(t)^\top Q(t), \quad X(t) = B(t)H(\mathcal{Z}(t), \tau(t)), \qquad (5.1)$$

is optimal for the problem (3.1)-(3.2), and $X(t)$ is the corresponding wealth with $X(0) = X_0$. This strategy replicates the claim $B(T)F(\mathcal{Z}(T), \lambda_J)$. If $\mathbf{E} U^+(F(\mathcal{Z}(T), \lambda_J)) < +\infty$ and either the uniqueness in Condition 3.1 holds for all $z$ or the law of $\mathcal{Z}_*(T)$ has a density, then the optimal strategy is unique.

**Proposition 5.1** Assume Conditions 3.1, 3.2, 3.3(ii) and 3.4. Then

$$H(x,t) = C_1 \left(\frac{x}{\lambda_J}\right)^\nu \exp\left\{\frac{1}{2}\nu(\nu-1)(T-t)\overline{R}\right\} + C_0$$

is the solution of (4.1) with $f(x) = F(x, \lambda_J)$ and the optimal strategy has the form

$$\pi(t)^\top = \nu B(t)(\widetilde{X}(t) - C_0)\widetilde{a}(t)^\top Q(t). \qquad (5.2)$$

It now follows that the solution of the problem (4.1), and hence the optimal strategy, can be written explicitly for the cases *(i)-(iv)* in Lemma 3.1. For the log utility case, $\nu = 1$ and, cf. above, the dependence of $H$ on $\overline{R}$ disappears; in fact $H$ is the steady state solution of (4.1). For cases *(ii)-(iv)* we require also that $R$ *be non-random* although we relax this for cases *(ii)* and *(iii)* in the next Corollary. For case *(iv)*, $\pi = \sum_j \pi^j$ where the $\pi^j$ are expressions of the form (5.2) with corresponding $\widetilde{X}^j$. In fact $\widetilde{X} = \sum_j \widetilde{X}^j$ and this decomposition depends on $\lambda_J$. For the case *(v)*, $H$ can be written in terms of the normal cumulative distribution function, so again the optimal strategy can be solved explicitly provided $R$ is non-random.

**Corollary 5.1** *Assume Conditions 3.1 and 3.4, and assume Condition 3.2 under the conditional probability given $R$. Then the optimal strategy is*

$$\bar{\pi}(t)^\top = \nu B(t)(\widetilde{X}(t) - C_0)\tilde{a}(t)^\top Q(t). \tag{5.3}$$

*where the normalized wealth $\widetilde{X}(t) = \widetilde{X}(t, \bar{\pi}(\cdot))$ is given by*

$$\begin{aligned} d\widetilde{X}(t) &= \nu(\widetilde{X}(t) - C_0)\tilde{a}(t)^\top Q(t)\sigma(t)dw_*(t), \\ \widetilde{X}(0) &= X_0. \end{aligned} \tag{5.4}$$

This result is a generalization of the case of "totally unhedgable" coefficients, cf. Karatzas and Shreve (1998), Chapter 6, Example 7.4. We see that the result holds for a larger class of utility functions than just the power utility functions, and the independence of the parameters and the Brownian motion can be relaxed considerably (when utility is only derived from terminal consumption). In fact, the corollary applies to cases *(i)-(iii)* in Lemma 3.1.

Here is another example where our theory applies, i.e. the optimal investment strategy is myopic, although this is not apparent from the corresponding Bellman equation. Assume that $\tilde{a}(t) \equiv \tilde{a}$ is constant in $t$ and independent of $w(\cdot)$. Let $\bar{\sigma}_i$, $i = 1, 2$, be given random matrices in $\mathbf{R}^{n \times n}$ which are independent of $w(\cdot)$, and let $\varepsilon > 0$ be fixed. Let $\tau'$ be any Markov time with respect to $\mathcal{F}_t$. Assume that

$$\sigma(t) \triangleq \begin{cases} \bar{\sigma}_1, & \text{if } t \notin [\tau, \tau + \varepsilon] \\ \bar{\sigma}_2, & \text{if } t \in [\tau, \tau + \varepsilon] \end{cases}, \quad \text{where} \quad \tau \triangleq \tau' \wedge (T - \varepsilon).$$

Then $R$ does not depend on $w(\cdot)$ (though $\mu(\cdot)$ does). For suitable $U$ we may apply Corollary 5.1 to obtain a strategy depending only on current wealth. If $\tilde{a}$ and $\bar{\sigma}_i$ are non-random, then $R$ is non-random, the market is still incomplete and according to Theorem 5.1 the strategy is myopic.

# 6 Optimal strategies for $m < n$

We turn now to the case of limited diversification. An individual investor may feel that she can only reasonably keep track of a limited number of stocks, or with finite capital wishes not to spread the investments too thinly, hence decides to hold at most $m$ stocks in her portfolio at any one time. A function $I \in \mathcal{I}_m$, defined in Section 2, will specify which $m$ stocks she holds at any time. If we restrict $\mathcal{I}_m$ to consist of constant functions only, then this amounts to choosing the $m$ "best" stocks initially and then trading in the market consisting of these $m$ stocks only. On the other hand, we may take $\mathcal{I}_m$ to consist of $\mu$-adapted processes taking values in $\mathcal{M}_m$. This form of $\mathcal{I}_m$ is not unreasonable. In fact, the rational investor, when choosing her portfolio, will want to maximize potential return while minimizing risk. As these factors depend only on the coefficient processes, $\mu$, it is reasonable to assume that $I$ is $\mu$-adapted. In this case, if the parameters $\mu$ are non-random, then the functions $I \in \mathcal{I}_m$ will be non-random but possibly time varying.

Recall that we write $\mu_I(t)$ for $(r(t), a_I(t), \sigma(t))$. We strengthen Conditions 3.2, 3.4 somewhat.

**Condition 6.1** *For all $I \in \mathcal{I}_m$, there exist $\lambda_{J_I} \triangleq \lambda_I \in \Lambda$, $C > 0$ and $c_0 \in (0, 1/(2J_I))$, such that $F(\cdot, \lambda_I)$ is piecewise continuous on $(0, \infty)$, $F(\mathcal{Z}_I(T), \lambda_I)$ is $\mathbf{P}_{*I}$-integrable, and*

$$\begin{cases} \mathbf{E}_{*I}\{F(\mathcal{Z}_I(T), \lambda_I)\} = X_0, \\ |F(z, \lambda_I)| \leq C z^{c_0 \log z} \quad \forall z > 0. \end{cases} \tag{6.1}$$

**Condition 6.2** *For all $I \in \mathcal{I}_m$,*

(i) *there exist $\lambda_{R_I} \triangleq \lambda_I \in \Lambda$ a.s., $C > 0$ and $c_{R_I} \triangleq c_I \in (0, 1/(2R_I))$ a.s., such that $F(\cdot, \lambda_I)$ is a.s. piecewise continuous on $(0, \infty)$, $F(\mathcal{Z}_I(T), \lambda_I)$ is a.s. $\mathbf{P}_{*I}(\cdot \mid R_I)$-integrable, and*

$$\begin{cases} \mathbf{E}_{*I}\{F(\mathcal{Z}_I(T), \lambda_I) \mid R_I\} = X_0 \text{ a.s.}, \\ |F(z, \lambda_I)| \leq C z^{c_I \log z} \quad \forall z > 0 \text{ a.s.} \end{cases} \tag{6.2}$$

(ii) $F(x, \lambda) = C_1(\lambda) x^\nu + C_0$, *where $C_1(\lambda) \neq 0$, $C_0$ and $\nu \neq 0$ are constants.*

We say that $\hat{I}$ *dominates* $I$ if $R_{\hat{I}} \geq R_I$, a.s. and $\mathbf{P}\{R_{\hat{I}} > R_I\} > 0$.

**Theorem 6.1** *Let $\hat{I} \in \mathcal{I}_m$. Assume Conditions 3.1, 6.1 and either*

(i) $U(x) = \log x$

or

(ii) $U^+(x) \leq \text{const}\,(|x|+1)$, $\mu_I$ and $w$ are independent for all $I \in \mathcal{I}_m$, and the random variable $R_{\hat{I}}$ is constant.

Then the strategy $\pi_{\hat{I}}(\cdot)$, defined in Theorem 5.1 (with $a(\cdot) = a_{\hat{I}}(\cdot)$) belongs to the class $\overline{\Sigma}(m)$ and

$$\mathbf{E}U(\widetilde{X}(T, \pi_{\hat{I}}(\cdot))) > \mathbf{E}U(\widetilde{X}(T, \pi(\cdot)))$$

for any strategy $\pi \in \overline{\Sigma}(m)$ such that $\hat{I}$ dominates the corresponding $I$.

Observe that $\mu_I$ amd $w$ are independent, in particular, if $\mathcal{I}_m$ contains only $\mu$-adapted processes, and $\mu(\cdot)$ and $w(\cdot)$ are independent.

**Corollary 6.1** *In the above theorem, the assumption in (ii) that $R_{\hat{I}}$ be non-random can be replaced by Condition 6.2 and then Condition 6.1 can be dropped.*

**Corollary 6.2** *Assume the hypotheses of either Theorem 6.1 or Corollary 6.1. If there exists $\hat{I} \in \mathcal{I}_m$ such that $R_{\hat{I}} = \max_{I \in \mathcal{I}_m} R_I$ a.s., then $\pi_{\hat{I}}(\cdot)$ is optimal for the problem (3.1)-(3.2). It is unique if $\hat{I}$ is.*

An interesting consequence is that the optimal $\hat{I} \in \mathcal{I}_m$ does not depend on $U(\cdot)$ or $\widehat{D}$ - just choose the $m$ stocks that provide a.s. the largest (in the $L_2$ sense) market price of risk. Because of the almost sure maximization requirement, this cannot always be done. However if $\mathcal{I}_m$ contains at least the $\mu$-adapted functions, then

$$\hat{I}(t) \in \arg\max_{M \in \mathcal{M}_m} \widetilde{a}(t)^\top Q_M(t) \widetilde{a}(t)$$

where $Q_M \triangleq Q_I$ with $I(\cdot) \equiv M$.

We can also apply our theory to a zero-coupon bond market based on a generalization of the Vasicek interest rate model. In the Vasicek model, the market price of risk is constant, cf. Lamberton and Lapeyre (1996), section 6.2.1. Our $\theta$ is their $-q$. We can generalize to $\theta$ a non-random function of $t$, so $R$ is non-random. Given a progressively measurable $u(\cdot): [0,T] \times \Omega \to [0,T]$ such that $u(t) \geq t$, a.s., we can construct the rolling bond $P(t, u(t))$ which expires at time $u(t)$. If $u(t) \equiv t$, this gives the usual bond-based construction of $B$. If $u(t) = [t] + 1$, this consists of a sequence of 1-year bonds rolled over at expiry. In a market consisting of the bank account $B$ and a finite number of such bonds and with a utility which satisfies Conditions 3.1, 3.2 as well as $U^+(x) \leq \text{const}\,(|x|+1)$, we can deduce that optimal strategies exist and one of them requires only one bond (with non-zero volatility) to be held in the optimal portfolio. As there is only one driving Brownian motion, this result is not surprising.

# 7 Appendix: Proofs

*Proof of Lemma 4.1.* By assumption, $\overline{R}$ and $J = T\overline{R}$ are non-random. Consider the following Cauchy problem:

$$\begin{cases} \frac{\partial V}{\partial t}(y,t) + \frac{\overline{R}}{2}\frac{\partial^2 V}{\partial y^2}(y,t) = 0, \\ V(y,T) = f(e^{y-J/2}), \quad y \in \mathbf{R}. \end{cases} \tag{7.1}$$

Let

$$p(y,t) \triangleq \frac{1}{\sqrt{2\pi(T-t)\overline{R}}} \exp\left(\frac{-y^2}{2\overline{R}(T-t)}\right). \tag{7.2}$$

It is the the fundamental solution of (7.1). By assumption, we have that

$$|f(e^y)| \leq Ce^{c_0 y^2} \quad \forall y \in \mathbf{R}, \quad c_0 < \frac{1}{2J} \leq \frac{1}{2\overline{R}(T-t)} \quad \forall t \in [0,T).$$

Then the integral

$$V(y,t) = \int_{-\infty}^{+\infty} p(y-y_0,t)f(e^{y_0-J/2})dy_0 \tag{7.3}$$

converges, and $V(y,t) \in C^{2,1}(\mathbf{R} \times (0,T))$ is a solution of the problem (7.1) such that $V(y,t) \to f(e^{y-J/2})$ a.s as $t \to T-$. Then $H(x,t) \triangleq V(\ln x + \overline{R}t/2, t)$ is a solution of (4.1) in the desired sense.

Note that $\tau(T) = T$ and $d\tau(t)/dt = |\theta(t)|^2/\overline{R}$. Set $\overline{H}(x,t) \triangleq H(x,\tau(t))$. Then

$$\begin{cases} \frac{\partial \overline{H}}{\partial t}(x,t) + \frac{|\theta(t)|^2}{2}x^2\frac{\partial^2 \overline{H}}{\partial x^2}(x,t) = 0, \\ \overline{H}(x,T) = f(x). \end{cases} \tag{7.4}$$

We define

$$\pi(t)^\top \triangleq B(t)\frac{\partial \overline{H}}{\partial x}(\mathcal{Z}(t),t)\mathcal{Z}(t)\widetilde{a}(t)^\top Q(t), \quad X(t) \triangleq B(t)\overline{H}(\mathcal{Z}(t),t).$$

Let us show that $X$ is the wealth, i.e. $X(t) = X(t, \pi(\cdot))$. Applying Itô's formula to $B(t)\overline{H}(\mathcal{Z}(t),t)$ and using (2.8), gives

$$dX(t) = r(t)X(t)\,dt + B(t)\frac{\partial \overline{H}}{\partial x}(\mathcal{Z}(t),t)\mathcal{Z}(t)\theta(t)^\top\,dw_*(t),$$

which is equivalent to (2.11) under our definition of $\pi$; hence $X$ is the wealth corresponding to $\pi$. If $\pi$ is admissible, then it replicates the claim $B(T)f(\mathcal{Z}(T))$, with initial wealth $X(0)$.

The integrability of $\pi$, cf Definition 2.2, follows if we take $T_k = T - 1/k$ and observe that $\frac{\partial \overline{H}}{\partial x}(\mathcal{Z}(t),t)$ is bounded pathwise for $t \leq T_k$ since $\frac{\partial \overline{H}}{\partial x}(x,t)$ is continuous on $(0,\infty) \times (0,T)$. The second requirement of Definition 2.2 follows from the continuity of $X(t) = B(t)\overline{H}(\mathcal{Z}(t),t)$, and the last follows from (4.3). So $\pi$ is admissible.

It remains to establish (4.3). Let

$$f_k(x) = \begin{cases} f(x) & \text{if } |f(x)| \leq k, \\ 0 & \text{otherwise.} \end{cases}$$

Let $\overline{H}_k(x,t)$ be the corresponding solution of (7.4) and let $X_k(t)$ be the corresponding wealth. Then $\widetilde{X}_k(t) = \overline{H}_k(\mathcal{Z}(t),t)$. Since $\lim_{t\uparrow T}\overline{H}_k(\mathcal{Z}(t),t) = \overline{H}_k(\mathcal{Z}(T),T) = f_k(\mathcal{Z}(T))$, a.s., and since $|\overline{H}_k(x,t)| \leq k$, then by dominated convergence $\lim_{t\uparrow T}\mathbf{E}_*\widetilde{X}_k(t) = \mathbf{E}_*f_k(\mathcal{Z}(T))$. Since $\widetilde{X}_k$ is a martingale for $t < T$, we then have $X_k(0) = \mathbf{E}_*f_k(\mathcal{Z}(T))$.

Further,

$$|\overline{H}(1,0) - \overline{H}_k(1,0)| = \int_{\{z:|f(e^{z-J/2})|>k\}} p(-z,0)|f(e^{z-J/2})|dz \to 0$$

since $\overline{H}(1,0) = \int p(-z,0)|f(e^{z-J/2})|dz < \infty$, cf. (7.2)-(7.3). Hence $X(0) = \overline{H}(1,0) = \lim_k \overline{H}_k(1,0) = \lim_k X_k(0) = \lim_k \mathbf{E}_*f(\mathcal{Z}(T)) = \mathbf{E}_*f(\mathcal{Z}(T))$. This completes the proof of Lemma 4.1. □

Let us establish a similar result in case Condition 3.3(i) holds.

**Lemma 7.1** *If Condition 3.3(i) is satisfied then $\widehat{\xi} \triangleq F(\mathcal{Z}(T), \lambda_J) = \mathcal{Z}(T)X_0$, $\lambda_J = X_0^{-1}$,*

$\mathbf{E}U(\widehat{\xi}) = J/2 + \log X_0$, and $B(T)\widehat{\xi}$ is replicated by $\pi$ where

$$\pi(t) = B(t)X_0 \mathcal{Z}(t)\widetilde{a}(t)^\top Q(t).$$

*Proof.* $F(z,\lambda) = z/\lambda$. From (3.4) we have $\lambda_J = X_0^{-1}$. Now $\mathbf{E}U(\widehat{\xi}) = \mathbf{E}\log[X_0 \mathcal{Z}(T)] = J/2 + \log X_0$. Let $\widetilde{X}(t) \triangleq \mathcal{Z}(t)X_0$. Then

$$\widetilde{X}(t) = X_0 + \int_0^T X_0 \mathcal{Z}(s)\theta(s)^\top \, dw_*(s)$$

and this is the normalized wealth corresponding to $\pi$ as given, cf (2.12). □

Let $U^-(x) \triangleq \max(0, -U(x))$; define a set of claims by

$$\Psi \triangleq \{\xi : \xi \text{ is } \mathcal{F}_T\text{-measurable}, \ \mathbf{E}_*|\xi| < +\infty, \ \mathbf{E}U^-(\xi) < \infty\},$$

and define $\mathcal{J}_i : \Psi \to \mathbf{R}$, $i = 0, 1$ by

$$\mathcal{J}_0(\xi) \triangleq \mathbf{E}U(\xi), \quad \mathcal{J}_1(\xi) \triangleq \mathbf{E}_*\xi - X_0.$$

Let us now define the claims *attainable in* $\widehat{D}$,

$$\Psi_0 \triangleq \{\xi \in \Psi : \xi = \widetilde{X}(T, \pi(\cdot)) \in \widehat{D} \text{ a.s.}, \ \pi(\cdot) \in \overline{\Sigma}\}.$$

Consider the problem

$$\text{Maximize} \quad \mathcal{J}_0(\xi) \quad \text{over} \quad \xi \in \Psi_0. \tag{7.5}$$

**Proposition 7.1** *Assume Conditions 3.1, 3.2 and 3.3. The optimization problem (7.5) has solution $\widehat{\xi} = F(\mathcal{Z}(T), \lambda_J)$.*

*Proof.* From Condition 3.2 it follows that $\mathbf{E}_*|\widehat{\xi}| < \infty$ and $\mathbf{E}_*\widehat{\xi} = X_0$. Let us show that $\mathbf{E}U^-(\widehat{\xi}) < \infty$. For $k = 1, 2, ...$, introduce the random events

$$\Omega^{(k)} \triangleq \{-k \leq U(\widehat{\xi}) \leq 0\},$$

along with their indicator functions $\chi^{(k)}$, respectively. The number $\widehat{\xi}$ achieves the unique maximum of the function $\mathcal{Z}(T)U(\xi) - \lambda_J \xi$ over the set $\widehat{D}$, and $X_0 \in \widehat{D}$. Hence for all $k = 1, 2, ...$,

$$\mathbf{E}\chi^{(k)}\left(U(\widehat{\xi}) - \mathcal{Z}_*(T)\lambda_J\widehat{\xi}\right) \geq \mathbf{E}\chi^{(k)}\left(U(X_0) - \mathcal{Z}_*(T)\lambda_J X_0\right) > -|U(X_0)| - |\lambda_J X_0|.$$

Since $\mathbf{E}\mathcal{Z}_*(T)|\widehat{\xi}| = \mathbf{E}_*|\widehat{\xi}| < +\infty$, then $\mathbf{E}U^-(\widehat{\xi}) < \infty$, hence $\widehat{\xi} \in \Psi$.

Let $L(\xi, \lambda) \triangleq \mathcal{J}_0(\xi) - \lambda \mathcal{J}_1(\xi)$, where $\xi \in \Psi$ and $\lambda \in \mathbf{R}$. We have

$$L(\xi, \lambda) = \mathbf{E}\left(U(\xi) - \lambda \mathcal{Z}_*(T)\xi\right) + \lambda X_0. \tag{7.6}$$

Let

$$\eta(\lambda) \triangleq F(\mathcal{Z}(T), \lambda) = F(\mathcal{Z}_*(T)^{-1}, \lambda). \tag{7.7}$$

From Condition 3.1 it follows that for any $\omega \in \Omega$, the random number $\eta(\lambda)$ provides the maximum in the set $\widehat{D}$ for the function under the expectation in (7.6).

Lemmas 4.1 and 7.1 imply the attainability of $\widehat{\xi}$, so $\widehat{\xi} \in \Psi_0$. Furthermore,

$$L(\xi, \lambda_J) \leq L(\widehat{\xi}, \lambda_J) \quad \forall \xi \in \Psi. \tag{7.8}$$

Let $\xi \in \Psi_0$ be arbitrary. We have that $\mathcal{J}_1(\xi) = 0$ and $\mathcal{J}_1(\widehat{\xi}) = 0$, then

$$\mathcal{J}_0(\xi) - \mathcal{J}_0(\widehat{\xi}) = \mathcal{J}_0(\xi) + \lambda_J \mathcal{J}_1(\xi) - \mathcal{J}_0(\widehat{\xi}) - \lambda_J \mathcal{J}_1(\widehat{\xi}) = L(\xi, \lambda_J) - L(\widehat{\xi}, \lambda_J) \leq 0.$$

Hence $\widehat{\xi}$ is an optimal solution of the problem (7.5). □

The following Lemma will prepare us for the proof of Theorem 5.1.

**Lemma 7.2** *Assume Conditions 3.1, 3.2, 3.3. Let $\widehat{\xi} \triangleq F(\mathcal{Z}(T), \lambda_J)$. Then*

(i) $\mathbf{E}U^-(\widehat{\xi}) < \infty$, $\widehat{\xi} \in \widehat{D}$ *a.s.*

(ii) *If $\mathbf{E}U^+(\widehat{\xi}) < +\infty$ and either the uniqueness in Condition 3.1 holds for all $z$ or the law of $\mathcal{Z}(T)$ has a density, then $\widehat{\xi}$ is unique (i.e. even for different $\lambda_J$, the corresponding $\widehat{\xi}$ agree a.s.).*

*Proof.* Part (i) is seen to hold from the proof of Proposition 7.1.

To show (ii), note that if $\mathbf{E}U^+(\widehat{\xi}) < +\infty$ then $L(\widehat{\xi}, \lambda_J) < +\infty$. Let $\xi' \in \Psi_0$ be an optimal solution of the problem (7.5). Let $\lambda_J$ be any number such that (3.4) holds. It is

easy to see that
$$L(\xi', \lambda_J) = \mathcal{J}_0(\xi') \geq \mathcal{J}_0(\widehat{\xi}) = L(\widehat{\xi}, \lambda_J).$$

From Condition 3.1 it follows that $\widehat{\xi} = \eta(\lambda_J)$ provides the maximum in the set $\widehat{D}$ of the function under the expectation in (7.6) with $\lambda = \lambda_J$. Hence both $\xi'$ and $\widehat{\xi}$ maximize $L(\cdot, \lambda_J)$. It follows that $\xi'$ must also maximize the function under the expectation in (7.6) a.s., and hence $\xi' = \eta(\lambda_J)$ a.s. from the uniqueness assertion in Condition 3.1. Thus (ii) is satisfied. This completes the proof of of Lemma 7.2. $\square$

*Proof of Theorem 5.1.* If $R = 0$ then $\mathcal{Z}_*(T) = 1$, hence the optimal claim $\widehat{\xi}$ is non-random. The only strategy which replicates the non-random claim is the trivial risk-free strategy, and (i) follows. It follows from Lemmas 7.1 and 4.1, that $\pi$ is admissible and replicates $B(T)\widehat{\xi}$. The optimality follows from Proposition 7.1.

Let us show the uniqueness of the optimal strategy. By Lemma 7.2 (ii), $\widehat{\xi}$ is the unique solution of the auxiliary problem (7.5). Hence if $\pi$ and $\pi'$ are two optimal strategies, they must both lead to the same wealth at time $T$. If we set $Y(t) = X(t, \pi'(\cdot)) - X(t, \pi(\cdot))$, then from (2.11) we obtain

$$dY(t) = r(t)Y(t)\,dt + (\pi'(t) - \pi(t))^\top \sigma(t)\,dw_*(t),\ Y(T) = 0.$$

Hence, given $\mu$, $(Y(t), (\pi'(t) - \pi(t))^\top \sigma(t))$ is a solution of the corresponding backward stochastic differential equation, $dY(t) = r(t)Y(t)\,dt + y(t)\,dw_*(t)$, $Y(T) = 0$. The theory of such equations, cf Yong and Zhou (1999), Theorem 2.2, p. 349, implies that the equation has a unique solution, $(Y, y) \equiv (0, 0)$. It follows that $\pi' = \pi$ a.e. a.s. $\square$

*Proof of Proposition 5.1.* Direct verification shows that $H$ is a solution of (4.1) and $\frac{\partial H(x,t)}{\partial x}x = \nu[H(x,t) - C_0]$. The result follows. $\square$

*Proof of Corollary 5.1.* Let $\overline{\Sigma}^R$ be the enlargement of $\overline{\Sigma}$ produced by replacing the filtration $\mathcal{F}_t$ by $\mathcal{F}_t^R$ generated by $\mathcal{F}_t$ and $R$ in the definition. With $\mathbf{P}(\cdot)$ replaced by $\mathbf{P}(\cdot|\mathcal{F}_0^R)$, we may apply Proposition 5.1 to obtain the optimal $\pi^R$ in feedback form (5.2) with $\widetilde{X}(t)$ replaced by $\widetilde{X}^R(t) = H(\mathcal{Z}(t), \tau(t))$. Note that $\lambda = \lambda_R$ now depends on $R$ and satisfies $X_0 = \mathbf{E}_*\{F(\mathcal{Z}(t), \lambda_R)|R\}$. Since $\widetilde{X}^R(t) = \widetilde{X}(t, \pi^R(\cdot))$ satisfies (5.4), then

$$\widetilde{X}(t, \pi^R(\cdot)) - C_0 = (X_0 - C_0)\exp\left\{\nu\int_0^t \theta(s)^\top dw_*(s) - \frac{\nu^2}{2}\int_0^t \|\theta(s)\|^2 ds\right\},$$

and hence $\widetilde{X}^R$ is $\mathcal{F}_t$ adapted and is *the* solution, $\widetilde{X}$, of (5.4).

Now set $\bar{\pi}(t)^\top = \mathbf{E}_*\{\pi^R(t)|\mathcal{F}_t\}^\top = \nu B(t)\left(\widetilde{X}(t) - C_0\right)\widetilde{a}(t)^\top Q(t)$. Then $\bar{\pi}$ is as defined in (5.3). Moreover $\bar{\pi}$ lies in the smaller control set $\overline{\Sigma}$, and

$$\begin{aligned}
\mathbf{E}U(\widetilde{X}(T,\bar{\pi}(\cdot)) &= \mathbf{E}_*\{\mathcal{Z}(T)U(\widetilde{X}(T,\bar{\pi}(\cdot)))\} = \mathbf{E}_*\{\mathcal{Z}(T)U(\widetilde{X}(T))\} \\
&= \mathbf{E}_*\{\mathcal{Z}(T)U(\widetilde{X}(T,\pi^R(\cdot)))\} = \mathbf{E}_*\mathbf{E}_*\{\mathcal{Z}(T)U(\widetilde{X}(T,\pi^R(\cdot)))|R\} \\
&\geq \mathbf{E}_*\mathbf{E}_*\{\mathcal{Z}(T)U(\widetilde{X}(T,\pi(\cdot)))|R\}
\end{aligned}$$

for any $\pi \in \overline{\Sigma}^R$. Hence $\bar{\pi}$ is optimal in this class, hence optimal for the original problem.
$\square$

*Proof of Theorem 6.1.* For the given $\hat{I}$, let us introduce a new market. Consider the auxiliary market defined by (2.1)-(2.2) with $a(\cdot)$ replaced by $a_{\hat{I}}(\cdot)$; we shall call it the $\hat{I}$–*market*. For $\pi \in \overline{\Sigma}(m)$ with corresponding $I = \hat{I}$, we have $\pi(t)^\top = \pi(t)^\top P_{\hat{I}}(t)$, and from (2.14) we have $P_{\hat{I}}(t)\widetilde{a}(t) = P_{\hat{I}}(t)\widetilde{a}_{\hat{I}}(t)$, so

$$\pi(t)^\top P_{\hat{I}}(t)[\sigma(t)\,dw(t) + \widetilde{a}_{\hat{I}}(t)\,dt] = \pi(t)^\top P_{\hat{I}}(t)[\sigma(t)\,dw(t) + \widetilde{a}(t)\,dt].$$

It follows from (2.12) that the wealth which is obtained with the strategy $\pi(\cdot)$ is the same for both markets – for the original market and for the $\hat{I}$–*market* – even though $w_*$ is different in the two markets.

Note that the assumptions of the Theorem suffice for uniqueness of the optimal strategy, cf. Theorem 5.1, because under *(i)* the minimum is unique, and under *(ii)*, $\mathcal{Z}_{*\hat{I}}(T)$ under $\mathbf{P}$ is conditionally log-normal given $\mu$ with parameters depending only on the constant $R_{\hat{I}}$, hence is unconditionally log-normal, hence has a density. Then we can apply Theorem 5.1 to the $\hat{I}$–*market* to obtain the unique optimal strategy $\pi_{\hat{I}} \in \overline{\Sigma}$. We show first that $\pi_{\hat{I}} \in \overline{\Sigma}(m)$.

$$\begin{aligned}
\pi_{\hat{I}}(t) &= B(t)\frac{\partial H}{\partial x}\Big(\mathcal{Z}_{\hat{I}}(t), \tau_{\hat{I}}(t)\Big)\mathcal{Z}_{\hat{I}}(t)Q(t)\widetilde{a}_{\hat{I}}(t) \\
&= B(t)\frac{\partial H}{\partial x}\Big((\mathcal{Z}_{\hat{I}}(t), \tau_{\hat{I}}(t)\Big)\mathcal{Z}_{\hat{I}}(t)Q(t)V(t)Q_{\hat{I}}(t)P_{\hat{I}}(t)\widetilde{a}(t) \qquad (7.9) \\
&= B(t)\frac{\partial H}{\partial x}\Big(\mathcal{Z}_{\hat{I}}(t), \tau_{\hat{I}}(t)\Big)\mathcal{Z}_{\hat{I}}(t)Q_{\hat{I}}(t)P_{\hat{I}}(t)\widetilde{a}(t).
\end{aligned}$$

We have used that $Q(t) = V(t)^{-1}$. Since $Q_{\hat{I}}(t)$ maps $L_{\hat{I}}(t)$ into $L_{\hat{I}}(t)$, then $\pi_{\hat{I}}(t) \in L_{\hat{I}}(t)$ for all $t$, so $\pi_{\hat{I}}(\cdot) \in \overline{\Sigma}(m)$, i.e. $P_{\hat{I}}(t)\pi_{\hat{I}}(t) \equiv \pi_{\hat{I}}(t)$. Let $\widetilde{X}_{\hat{I}}(t)$ be the corresponding normalized wealth. By Theorem 5.1, there exists $\lambda_{J_{\hat{I}}} \triangleq \lambda_{\hat{I}} \in \Lambda$ such that $\widetilde{X}_{\hat{I}}(T) = F(\mathcal{Z}_{\hat{I}}(T), \lambda_{\hat{I}})$.

Now assume $\hat{I}$ dominates $I$ and consider a new auxiliary market which we shall call the $I^+$–market: we assume that this market consists of the bond $B(t)$ and the stocks $S_1(t), ..., S_n(t), S_{n+1}(t)$, where the stock prices $S_1(t), ..., S_n(t)$ are defined by (2.1), replacing $a(\cdot)$ by $a_I(\cdot)$, and where $S_{n+1}(t)$ is defined by the equation

$$dS_{n+1}(t) = S_{n+1}(t)\left((r(t) + \alpha)dt + dw_{n+1}(t)\right), \tag{7.10}$$

with

$$\alpha \triangleq \sqrt{\frac{R_{\hat{I}} - R_I}{T}}, \tag{7.11}$$

and with $w_{n+1}(t)$ a scalar Wiener process independent of $(w(\cdot), \mu(\cdot))$. Of course the filtration for this market, $\{\mathcal{F}_t^+\}$, will be larger than $\{\mathcal{F}_t\}$ (it includes information on $w_{n+1}$ and $\alpha$). It is easy to see that the corresponding numbers $J_{I^+}$ and $R_{I^+}$ for the $I^+$–market are

$$R_{I^+} = R_I + \alpha^2 T = R_{\hat{I}}, \quad J_{I^+} = J_I + \alpha^2 T = J_{\hat{I}}. \tag{7.12}$$

It follows that if (ii) holds, then the distribution of $\mathcal{Z}_{I^+}(T)$ under $\mathbf{P}_{*I^+}$ and of $\mathcal{Z}_{\hat{I}}(T)$ under $\mathbf{P}_{*\hat{I}}$ are both log-normal with mean equal to 1 and variance of the log equal to $R_{I^+} = R_{\hat{I}}$, hence the same. Thus

$$\mathbf{E}_{*I^+}\{|F(\mathcal{Z}_{I^+}, \lambda_{\hat{I}})|\} = \mathbf{E}_{*\hat{I}}\{|F(\mathcal{Z}_{\hat{I}}, \lambda_{\hat{I}})|\} < +\infty,$$

and by Lemma 4.1

$$\mathbf{E}_{*I^+}\{F(\mathcal{Z}_{I^+}, \lambda_{\hat{I}})\} = \mathbf{E}_{*\hat{I}}\{F(\mathcal{Z}_{\hat{I}}, \lambda_{\hat{I}})\} = X_0.$$

Lemma 7.1 implies the same result if (i) holds. Then Conditions 3.1-3.3 are satisfied for the $I^+$–market. By Theorem 5.1 applied to the $I^+$–market, there exists a unique optimal strategy for the problem (3.1)-(3.2) for the class $\overline{\Sigma}^+$ for this market. Let $\widetilde{X}_{I^+}(t)$ be the corresponding normalized wealth for the $I^+$–market.

In case (i), Lemma 7.1 and (7.12), imply $\mathbf{E}U(\widetilde{X}_{I^+}(T)) = J_{\hat{I}}/2 + \log X_0 = \mathbf{E}U(\widetilde{X}_{\hat{I}}(T))$.

Similarly, in case *(ii)*, (7.12) implies

$$
\begin{aligned}
\mathbf{E}U(\widetilde{X}_{I^+}(T)) &= \mathbf{E}_{*I^+}\{\mathcal{Z}_{I^+}(T)U(F(\mathcal{Z}_{I^+}(T),\lambda_{\hat{I}}))\} \\
&= \mathbf{E}_{*\hat{I}}\{\mathcal{Z}_{\hat{I}}(T)U(F(\mathcal{Z}_{\hat{I}}(T),\lambda_{\hat{I}}))\} \\
&= \mathbf{E}U(\widetilde{X}_{\hat{I}}(T)).
\end{aligned}
\quad (7.13)
$$

Now consider an arbitrary strategy $\pi(\cdot) \in \overline{\Sigma}(m)$, with corresponding $I$ dominated by $\hat{I}$, as a strategy in the $I^+$–market, with the investment in stock $(n+1)$ equal zero identically. By Theorem 5.1 the unique optimal strategy in the $I^+$–market holds a non-zero multiple of $\alpha$ in stock $(n+1)$, hence $\pi(\cdot)$ is not optimal. Then by (7.13)

$$\mathbf{E}U(\widetilde{X}(T,\pi)) < \mathbf{E}U(\widetilde{X}_{I^+}(T)) = \mathbf{E}U(\widetilde{X}_{\hat{I}}(T)).$$

This completes the proof of Theorem 6.1. □

*Proof of Corollary 6.1.* For given $I \in \mathcal{I}_m$ let us enlarge the filtration $\{\mathcal{F}_t\}$ to $\{\mathcal{F}_t^{R_I}\} \triangleq \{\mathcal{F}_t^I\}$ with corresponding $\overline{\Sigma}^I$. Recall that $w$ is still Brownian motion because of the independence of $w$ and $\mu_I$. We will work with the conditional probability, $\mathbf{P}(\cdot|\mathcal{F}_0^I)$, rather than with $\mathbf{P}$. Uniqueness of the optimal strategy still holds under under this measure.

As in the proof of Corollary 5.1 applied to the $\hat{I}$-market, the strategy $\pi_{\hat{I}}(t) = \nu B(t)[\widetilde{X}(t) - C_0]Q(t)\widetilde{a}_{\hat{I}}(t) \in \overline{\Sigma}$ is optimal in $\overline{\Sigma}^{\hat{I}} \supset \overline{\Sigma}$. Moreover, as in (7.9), $\pi_{\hat{I}} \in \overline{\Sigma}(m)$.

Proceeding as in the proof of Theorem 6.1, we can define the $I^+$-market. Then $\mathcal{F}_t^+$ is generated by $(\mu_I(\cdot), \alpha, w(\cdot), w_{n+1}(\cdot))$, and $\mathcal{F}_t^{I^+}$ by $(\mu_I(\cdot), \alpha, w(\cdot), w_{n+1}(\cdot), R_{I^+})$. The corresponding set of policies is denoted by $\overline{\Sigma}^{I^+}$. Again $R_{I^+} = R_{\hat{I}}$ and the conditional distribution of $\mathcal{Z}_{I^+}$ under $\mathbf{P}_{*I^+}$ given $R_{I^+}$ is the same as the conditional distribution of $\mathcal{Z}_{\hat{I}}$ under $\mathbf{P}_{*\hat{I}}$ given $R_{\hat{I}}$. Then (7.13) becomes

$$
\begin{aligned}
\mathbf{E}\{U(\widetilde{X}_{I^+}(T)) \mid R_{I^+}\} &= \mathbf{E}_{*I^+}\{\mathcal{Z}_{I^+}(T)U(F(\mathcal{Z}_{I^+}(T),\lambda_{\hat{I}})) \mid R_{I^+}\} \\
&= \mathbf{E}_{*\hat{I}}\{\mathcal{Z}_{\hat{I}}(T)U(F(\mathcal{Z}_{\hat{I}}(T),\lambda_{\hat{I}})) \mid R_{\hat{I}}\} \\
&= \mathbf{E}\{U(\widetilde{X}_{\hat{I}}(T)) \mid R_{\hat{I}}\}.
\end{aligned}
$$

Taking expectations gives

$$\mathbf{E}\{U(\widetilde{X}_{I^+}(T))\} = \mathbf{E}\{U(\widetilde{X}_{\hat{I}}(T)).\} \quad (7.14)$$

Any $\pi \in \overline{\Sigma}(m)$, with corresponding $I$ dominated by $\hat{I}$, can be considered as as a (non-optimal) element of $\Sigma^{I^+}$ and so (7.14) implies

$$\mathbf{E}U(\widetilde{X}(T,\pi)) < \mathbf{E}U(\widetilde{X}_{I^+}(T)) = \mathbf{E}U(\widetilde{X}(T,\pi_{\hat{I}})). \qquad \Box$$

The proof of Corollary 6.2 is obvious.

# References

[1] Black, F., Scholes, M.: The pricing of options and corporate liabilities, J. Political Econ. **81**, 637-659 (1973).

[2] Karatzas, I., Shreve, S.E.: Brownian Motion and Stochastic Calculus, 2nd edn, Springer-Verlag, New York (1991).

[3] Karatzas, I., Shreve, S.E.: Methods of Mathematical Finance, Springer-Verlag, New York, (1998).

[4] Khanna, A., Kulldorff, M.: A generalization of the mutual fund theorem, Finance Stochast. **3**, 167-185 (1999).

[5] Laurent, J.P., Pham, H.: Dynamic programming and mean-variance hedging, Finance Stochast. **3**, 83-110 (1999).

[6] Lambertone, D., Lapeyre, B: Introduction to Stochastic Calculus Applied to Finance, Chapman & Hall, London (1996).

[7] Merton, R.: Lifetime portfolio selection under uncertainty: the continuous-time case, Rev. Econ. Statist. **51**, 247-257 (1969).

[8] Merton, R.: The theory of rational option pricing. Bell J. Econ. Manag. Sci. **4**, 141-183 (1973).

[9] Merton, R.: Continuous-Time Finance, Blackwell, Cambridge, Mass., (1990).

[10] Yong, J., Zhou, X. Y.: Stochastic Controls, Springer-Verlag, New York (1999).

[11] Zhou, Z.: An equilibrium analysis of hedging with liquidity constraints, speculations, and government price subsidy in a commodity market, J. Finance **53**, 1705-1736 (1998).